\documentclass[11pt]{amsart}

\usepackage{fullpage}
\usepackage{amsmath}
\usepackage{amssymb}
\usepackage{verbatim}
\usepackage[linktocpage]{hyperref}

\newtheorem{df}{Definition}[section]
\newtheorem{thm}[df]{Theorem}

\newtheorem{rem}[df]{Remark}

\newtheorem{lem}[df]{Lemma}

\newcommand{\pf}{\textit{Proof.} }
\newcommand{\eq}{eqnarray*}

\newcommand{\aK}{almost K\"{a}hler }
\newcommand{\aH}{almost Hermitian }

\begin{document}

\title{Lower Order Tensors in Non-K\"{a}hler Geometry and Non-K\"{a}hler Geometric Flow}
\date{}
\author{Song Dai}


\address{Center for Applied Mathematics of Tianjin University\\
Tianjin University\\
No.92 Weijinlu Nankai District\\
Tianjin\\
P.R.China 300072}

\email{song.dai@tju.edu.cn}

\begin{abstract}
In recent years, Streets and Tian introduced a series of curvature flows to study non-K\"{a}hler geometry. In this paper, we study how to construct the second order curvature flows in a uniform way, under some natural assumptions which hold in Streets and Tian's works. As a result, by classifying the lower order tensors, we classify the second order curvature flows in \aH, \aK and Hermitian geometries in certain sense. In particular, the Symplectic Curvature Flow is the unique way to generalize Ricci Flow on \aK manifolds.
\end{abstract}

\maketitle

\tableofcontents

\section{Introduction}
In 1982, Hamilton \cite{Ham} introduced the Ricci Flow $\frac{\partial}{\partial t}g=-2Ric$ on Riemannian manifolds and showed that compact $3$-manifold with positive Ricci curvature is the spherical space form. After that, many people studied the Ricci Flow intensively. In 2002, Perelman \cite{Pe1}\cite{Pe2}\cite{Pe3} did breakthrough that by using the Ricci Flow, he proved Thurston's Geometric Conjecture and as a corollary, Poincar\'{e} Conjecture. Since the method of curvature flow is so powerful, people tried to use the similar idea to study other geometric objects. In 1985, Cao \cite{Cao} initialized the study of the Ricci Flow on K\"{a}hler manifolds, which is the K\"{a}hler Ricci Flow. He showed that if we fix the complex structure $J$, the Ricci Flow preserves the K\"{a}hler structure.

To generalize the K\"{a}hler Ricci Flow to the non-K\"{a}hler case, Streets and Tian introduced a series of the second order curvature flows on this subject, including the Hermitian Curvature Flow \cite{ST her}, the Symplectic Curvature Flow \cite{ST symp}, the Almost Hermitian Curvature Flow \cite{ST symp}, the Pluriclosed Flow \cite{ST pluri}, the Pluriclosed Flow on generalized K\"{a}hler manifolds (or Generalized K\"{a}hler Ricci Flow) \cite{ST gkg}. Along this direction, other people constructed new flows. Vezzoni \cite{V} generalized the Hermitian Curvature Flow to almost Hermitian setting. The author \cite{Dai} unified the Symplectic Curvature Flow and the Pluriclosed Flow in almost Hermitian setting. For other relevant studies , one may refer to \cite{B}\cite{E}\cite{EFV}\cite{F}\cite{FV}\cite{Lau}\cite{LW}\cite{Po}\cite{Sm}\cite{Streets pgs}\cite{Streets pbg}\cite{Streets gkrf}\cite{ST reg pluri}\cite{SW}. Besides this series of works, there are other ways to define the second order curvature flows on non-K\"{a}hler geometry. Gill, Tosatti and Weinkove \cite{Gill}\cite{TW} studied the Chern Ricci Flow on Hermitian manifolds. Wang and L\^{e} \cite{LeWang} studied the Anti-Complexified Ricci Flow on \aK manifolds.

In this paper, we focus on three kinds of non-K\"{a}hler geometries: almost Hermitian geometry, \aK geometry and Hermitian geometry. We discuss how to construct the second order curvature flows in certain ``canonical" sense. Basically, we show that the geometric flows defined by Streets and Tian in above three geometries have some ``canonical" uniqueness in certain sense.
\\

Let $(g,J,\omega)$ be an almost Hermitian structure. Let $T=T(g,J,\omega)$ be a tensor defined from $(g,J,\omega)$. Notice that $(g,-J,-\omega)$ is also an almost Hermitian structure. We say that $T$ is of even type, if $T(g,J,\omega)=T(g,-J,-\omega)$, for example $Ric(X,Y)$. We say that $T$ is of odd type, if $T(g,J,\omega)=-T(g,-J,-\omega)$, for example $Ric(JX,Y)$. We say that $T$ scales as $r^{k}$, if $T(r^{k}g,J,r^{k}\omega)=T(g,J,\omega)$.

Let $T$ be a $2$-tensor. Consider the $J$ action on $T$ given by $J^{*}T(X,Y)=T(JX,JY)$. We denote
\begin{\eq}
T^{(1,1)}(X,Y)&:=&\frac{1}{2}(T(X,Y)+T(JX,JY)),\\
T^{(0,2)+(2,0)}(X,Y)&:=&\frac{1}{2}(T(X,Y)-T(JX,JY)).
\end{\eq}
Consider the transposition action on $T$ given by ${}^{t}T(X,Y)=T(Y,X)$. We denote
\begin{\eq}
T^{sym}(X,Y)&:=&\frac{1}{2}(T(X,Y)+T(Y,X)),\\
T^{skew}(X,Y)&:=&\frac{1}{2}(T(X,Y)-T(Y,X)).
\end{\eq}

For tensors $A,B$, we denote by $A*B$ any linear combination of $g$-operations on $A\otimes B$. By a $g$-operation we mean any composite of raising or lowering indices and taking traces. 

Let $(g_{t},J_{t},\omega_{t})$ be a family of almost Hermitian structures. For simplicity of writing we drop the subscript $t$ immediately and define
\begin{\eq}
\frac{\partial}{\partial t}g=h,\quad\frac{\partial}{\partial t}J=K,\quad \frac{\partial}{\partial t}\omega=\eta.
\end{\eq}
Notice that every two of $(h,K,\eta)$ determines the third one. We may focus on the pair $(h,K)$.
\\

We construct the geometric flows under the following natural assumption.
\\
$\mathbf{Assumption}$ (A)\\
(1) $(h,K)$ is tensorial in $(g,J,\omega)$ and is a differential operator with respect to $(g,J,\omega)$.\\
(2) $(h,K)$ is a second order system with respect to $(g,J)$.\\
(3) $h$ is of even type, and $K$ is of odd type.\\
(4) $h$ scales as $r^{0}$ and $K$ scales as $r^{-2}$.\\
(5) Modulo gauge transformation, the symbol of the system $(h,K)$ is identity, where the gauge transformation is generated by a vector field defined from $(g,J,\omega,\overline{D})$, where $\overline{D}$ is a fixed linear connection.

For Assumption (A), (1) is clearly natural. Since we only consider the second order flows, we require (2). For (3), we hope that the flows should be invariant if we replace $(g,J)$ by $(g,-J)$. We require (4) since we hope the flows should be invariant under parabolic rescaling, i.e., $(g_t,J_t)\mapsto (r^2g_{r^2t},J_{r^2t})$. The only essential assumption is (5), which ensures short-time existence by De Turck's trick. This is indeed the case for the Ricci Flow. So we regard Assumption (A) as a natural and basic assumption to generalize the K\"{a}hler Ricci Flow to non-K\"{a}hler case.

Most non-K\"{a}hler geometric flows I mentioned above satisfy Assumption (A), including Streets and Tian's series of works, the Hermitian Curvature Flow \cite{ST her}, the Symplectic Curvature Flow \cite{ST symp}, the Almost Hermitian Curvature Flow \cite{ST symp}, the Pluriclosed Flow \cite{ST pluri}, and the flow introduced by Vezzoni \cite{V}, the flow introduced by the author \cite{Dai}, the Chern Ricci Flow studied by Gill, Tosatti and Weinkove \cite{Gill}\cite{TW}. In the Hermitian setting, the flows above satisfy condition (5') introduced below instead of (5). But we may take a gauge transformation using the Lee form to obtain condition (5), see Remark \ref{hergauge}.
\\

Under Assumption (A), we classify the second order curvature flows in the following non-K\"{a}hler geometries.

In almost Hermitian geometry, we have the following result.
\begin{thm}\label{ah}
In almost Hermitian geometry, under Assumption (A), the geometric flows are of the following form.
\begin{\eq}
\frac{\partial}{\partial t}g&=&-2Ric+aL_{\theta^{\sharp}}g+Q_{1}\\
\frac{\partial}{\partial t}J&=&\triangle J+\mathcal{N}+\mathcal{R}+aL_{\theta^{\sharp}}J+Q_{2},
\end{\eq}
where $a\in \mathbb{R}$, $Q_{1},Q_{2}$ are of the form $DJ*DJ$ and satisfy the algebraic, necessary conditions
\begin{\eq}
&&Q_{1}~\text{is symmetric},\\
&&Q_{2}~\text{is}~(0,2)+(2,0),\\
&&Q_{1}^{(0,2)+(2,0)}=Q_{2}^{\text{sym}}J.
\end{\eq}
\end{thm}
$(-2Ric, \triangle J+\mathcal{N}+\mathcal{R})$ is the second order term which guarantees the system will have a ``good" symbol, where $Ric$ is the Ricci curvature with respect to the Levi-Civita connection $D$, $\triangle$ is the rough Laplacian with respect to the Levi-Civita connection $D$,
\begin{\eq}
g(\mathcal{N} (X),Y)&=&g^{ab}g(D_{a}J(JX),D_{b}J(Y)), \\
g(\mathcal{R}(X),Y)&=&Ric(JX,Y)+Ric(X,JY), \\
\theta(X)&=&g^{ij}DJ(e_{i},Je_{j},X).
\end{\eq}
The gauge term is then $(aL_{\theta^{\sharp}}g, aL_{\theta^{\sharp}}J)$.

We notice that the Almost Hermitian Curvature Flow in \cite{ST symp} is in the family above.
\\

In \aK geometry, we need to require $d\eta=0$. One natural option is $\eta=P$, where $P$ is the Chern form up to a factor. It is known that in \aK setting the Chern connection is the unique Hermitian connection (see also the Appendix). Then we have the following result.
\begin{thm}\label{ak}
The only geometric flow in \aK geometry satisfying Assumption (A) and $\frac{\partial}{\partial t}\omega=P$ is the Symplectic Curvature Flow defined in \cite{ST symp}.
\begin{\eq}
\frac{\partial}{\partial t}\omega&=&P\\
\frac{\partial}{\partial t}g&=&\triangle J+\mathcal{N}+\mathcal{R}.
\end{\eq}
\end{thm}
Notice that in \cite{ST symp}, Streets and Tian mentioned that one may modify the Symplectic Curvature Flow by adding some first order terms. However, here we rule out this possibility. This theorem tells us that the Symplectic Curvature Flow is the unique way to generalize the Ricci Flow on \aK manifolds.
\\

In Hermitian geometry, we should require the integrability of $J$ is preserved. We would like to fix the complex structure $J$. In this case, the system is not gauge invariant, so we assume that:\\
(5') the symbol of $h$ is the identity with respect to $g$ (not modulo gauge transformation).

Notice that we may add the gauge term $L_{\theta^{\sharp}}J$ to obtain condition (5) from (5'), see remark \ref{hergauge}.
\begin{thm}\label{h}
In Hermitian geometry, if we assume $\frac{\partial}{\partial t} J=0$, under (1)(2)(3)(4) in Assumption (A) and (5') above, the geometric flows are of the following form
\begin{\eq}
\frac{\partial}{\partial t}\omega&=&S+\mathcal{Q}\\
\frac{\partial}{\partial t}J&=&0,
\end{\eq}
where $\mathcal{Q}=a_{1}B^{1}J+a_{2}B^{2}J+a_{3}(B^{5})^{(1,1)}J+a_{4}(B^{6})^{sym}J$.
\end{thm}
Here $S(X,Y)=g^{ij}g(\Omega(e_{i},Je_{j},X),Y)$, $\Omega$ is the curvature with respect to the Chern connection. And the definition of $B^{i}$ is given in section \ref{tensor}.
\\

To establish the above results, first we derive the algebraic necessary conditions and the analytic sufficient conditions for deforming \aH structures. Then to find the suitable tensors satisfying these conditions, we classify the lower order tensors in the corresponding geometry, more precisely, the first order tensors, the second order tensors and the gauge terms. Finally, by calculating the symbols, we find out the desired tensors.

We organize the paper as follows. In Section \ref{pre}, we recall some preliminaries in almost Hermitian geometry and fix some notations. In Section \ref{con}, we derive the algebraic necessary conditions and the analytic sufficient conditions for deforming almost Hermitian structures, which provide the main restrictions in constructing the curvature flows. In Section \ref{tensor}, we classify the lower order tensors in \aH geometry, \aK geometry and Hermitian geometry. In Section \ref{flow}, we calculate the symbols of the second order tensors and then classify the second order curvature flows in \aH, \aK and Hermitian geometries. In the Appendix, we recall two basic facts about Hermitian connections.
\\

\textbf{Acknowledgements:} The author wishes to express his gratitude to his advisor Gang Tian, for suggesting the author to study the problems in the non-K\"{a}hler geometric flow, especially the Symplectic Curvature Flow, and encouraging the author all the time and many helpful discussions. The author would also like to thank Jeffrey Streets for his helpful comments and suggestions.

\section{Preliminaries}\label{pre}
Let $M$ be a manifold, $J$ be a section of $End(TM)$. We say that $J$ is an almost complex structure, if $J^{2}=-1$. We say that an almost complex structure $J$ is integrable, if $J$ is induced by holomorphic coordinates. By the theorem of Newlander-Nirenberg \cite{NN}, $J$ is integrable if and only if $N_{J}=0$, where
\begin{eqnarray*}
N_{J}(X,Y)=[JX,JY]-[X,Y]-J[JX,Y]-J[X,JY]
\end{eqnarray*}
is the Nijenhuis tensor.

We say that $(g,J,\omega)$ is an almost Hermitian structure if the following conditions hold.\\
(1) $g$ is a Riemannian metric.\\
(2) $J$ is an almost complex structure.\\
(3) $(g,J)$ is compatible, i.e.,
$g(JX,JY)=g(X,Y)$.\\
(4) We define
$\omega(X,Y)=g(JX,Y).$

Moreover $(g,J, \omega)$ is Hermitian if $J$ is integrable, and $(g,J,\omega)$ is almost K\"{a}hler if $d\omega=0$. Furthermore $(g,J,\omega)$ is K\"{a}hler if $J$ is integrable and $d\omega=0$.

We fix some notations first.
\\
\textbf{Notations:}

Let $(g,J,\omega)$ be an almost Hermitian structure.
\\
(1) Let $D$ denote the Levi-Civita connection and $D$ is extended to tensor fields. For instance,
\begin{\eq}
DJ(X,Y)=(D_{X}J)Y=D_{X}(JY)-J(D_{X}Y).
\end{\eq}
(2) We implicitly identify $TM$ and $T^{*}M$ by using $g$, i.e., for instance,
\begin{\eq}
DJ(X,Y,Z)=g(DJ(X,Y),Z).
\end{\eq}

Notice that in this notation, $\omega=J$.
\\
(3) Usually, we use $i$ instead of $e_{i}=\frac{\partial}{\partial x^{i}}$ for short. We use orthonormal basis at one point, and often we assume it is normal. The same index means to take (real) trace with respect to $g$. For complex trace, we mean that $\omega^{ij}T_{ij}$, and it equals to $T(i,Ji)$, where $i$ goes over all the orthonormal basis.\\
\begin{rem}
Since $g$ is parallel with respect to our connection, the notation $(2)$ is safe if we do tensor calculation in a fixed geometric structure. But we should take care if we calculate the evolution equations.
\end{rem}
\begin{rem} We say that $T$ is $(1,1)$ ($(0,2)+(2,0)$), if $T^{(0,2)+(2,0)}=0$ ($T^{(1,1)}=0$). Then by our identification, in \aH setting,
\begin{\eq}
T\quad~\text{is}~ (1,1)&\Longleftrightarrow& TJ=JT,\\
T\quad~\text{is}~ (0,2)+(2,0)&\Longleftrightarrow& TJ=-JT.
\end{\eq}
\end{rem}
\pf
We only prove the second identity. By definition,
\begin{eqnarray*}
\langle(TJ+JT)X,Y\rangle=T(JX,Y)-T(X,JY)=2T^{(1,1)}(JX,Y).
\end{eqnarray*}
\qed

For contraction, we have the following lemma.
\begin{lem}
In the following cases, the tensors will vanish.\\
(a) Taking trace of a skew $2$-tensor.\\
(b) Taking complex trace of a symmetric $2$-tensor.\\
(c) Either taking trace or complex trace of a $(0,2)+(2,0)$ $2$-tensor.
\end{lem}
\pf
We only prove (c). Since $\{e_{i}\}$ is a orthonormal basis and $J$ is isometry, $\{Je_{i}\}$ is also a orthogonal basis. So to take trace, we may replace $i$ by $Ji$. By definition, we have
\begin{\eq}
T(i,i)&=&T(Ji,Ji)=T(JJi,i)=-T(i,i),\\
T(i,Ji)&=&T(Ji,i)=T(JJi,Ji)=-T(i,Ji).
\end{\eq}
\qed

We come back to preliminaries. Let $(g,J,\omega)$ be an almost Hermitian structure. Let $\nabla$ be the corresponding Chern connection, i.e.
\begin{\eq}
\nabla g=0,\quad \nabla J=0,\quad \text{Tor}^{\nabla}~ \text{is}~ (0,2)+(2,0)~ \text{for the first two variables},
\end{\eq}
where $\text{Tor}^{\nabla}(X,Y)=\nabla_{X}Y-\nabla_{Y}X-[X,Y]$ is the torsion of $\nabla$. From \cite{Ga} or Appendix, we see
\begin{\eq}
g(\nabla_{X}Y,Z)&=&g(D_{X}Y,Z)+\frac{1}{2}DJ(X,JY,Z)\\
&&+\frac{1}{4}(DJ(JY,Z,X)+DJ(JZ,X,Y)-DJ(Y,Z,JX)-DJ(Z,X,JY)).
\end{\eq}

And in \aK setting,
\begin{\eq}
g(\nabla_{X}Y,Z)&=&g(D_{X}Y,Z)+\frac{1}{2}DJ(X,JY,Z).
\end{\eq}

Let $\Omega$ denote the curvature of $\nabla$, i.e.
\begin{\eq}
\Omega(X,Y,Z,W)=g(\nabla_{X}\nabla_{Y}Z-\nabla_{Y}\nabla_{X}Z-\nabla_{[X,Y]}Z,W).
\end{\eq}

And we also define the Riemannian curvature $Rm$ in the same manner. Define
\begin{\eq}
P(X,Y)=\Omega(X,Y,i,Ji),\quad S(X,Y)=\Omega(i,Ji,X,Y),\quad Ric(X,Y)=Rm(i,X,Y,i).
\end{\eq}

Notice that $-\frac{1}{4\pi}P$ is Chern form, i.e., $[-\frac{1}{4\pi}P]=c_{1}$, where $c_{1}$ is the first Chern class of $T^{1,0}M$. And $S$ is a $(1,1)$ form.
\begin{rem}
Here $P$ differs a minus from $P$ in \cite{ST symp} since the definition of curvature differs a minus.
\end{rem}

Denote $\rho'(X,Y)=Rm(JX,Y,i,Ji)$, $s'=Rm(i,Ji,j,Jj)$. Let $\theta$ denote the Lee form, i.e., $\theta=DJ(i,Ji,\cdot)$.

Now we recall some basic identities in Riemannian geometry.\\
(1) Symmetries of $Rm$.\\
(2) $1$st and $2$nd Bianchi identity:
\begin{\eq}
&1st:&~ Rm(X,Y,Z,W)+Rm(Y,Z,X,W)+Rm(Z,X,Y,W)=0.\\
&2nd:&~ DRm(X,Y,Z,W,V)+DRm(Y,Z,X,W,V)+DRm(Z,X,Y,W,V)=0.
\end{\eq}
(3) Ricci identity: For example, let $T$ be a $2$-tensor, then
\begin{\eq}
D^{2}T(X,Y,Z,W)-D^{2}T(Y,X,Z,W)&=&(Rm(X,Y)T)(Z,W)\\
&=&-T(Rm(X,Y,Z),W)-T(Z,Rm(X,Y,W)).
\end{\eq}

Then we consider some basic identities in \aH geometry. For the following material, one may refer to \cite{Ga}.

Let $(g,J,\omega)$ be an almost Hermitian structure.
\begin{lem}\label{dj}
$DJ$ is skew and $(0,2)+(2,0)$ with respect to the last two components.
\begin{\eq}
DJ(X,Y,Z)&=&-DJ(X,Z,Y)\\
DJ(X,JY,JZ)&=&-DJ(X,Y,Z)
\end{\eq}
\end{lem}
\pf
By our notation, $DJ=D\omega$, so the first identity holds. For the second one, fix a point $p$, suppose $X,Y,Z$ are in a normal frame at $p$, then
\begin{\eq}
&0&=D_{X}g(Y,Z)=D_{X}(g(JY,JZ))=g(D_{X}J(Y),JZ)+g(JY,D_{X}J(Z))\\
&&=DJ(X,Y,JZ)+DJ(X,Z,JY)=DJ(X,Y,JZ)-DJ(X,JY,Z).
\end{\eq}

Then we finish the proof.
\qed
\begin{lem}\label{hc}
(Hermitian condition) Let $(g,J,\omega)$ be an \aH structure. Then
\begin{\eq}
N_{J}=0&\Longleftrightarrow& DJ(JX,Y,Z)-DJ(JY,X,Z)+DJ(X,JY,Z)-DJ(Y,JX,Z)=0\\
&\Longleftrightarrow& DJ(JX,JY,Z)=DJ(X,Y,Z)
\end{\eq}
\end{lem}
\pf
The first identity is from the definition of $N_{J}$ and $D_{X}Y-D_{Y}X=[X,Y]$. For the second identity, adding the following identities together,
\begin{\eq}
DJ(JX,Y,Z)-DJ(JY,X,Z)+DJ(X,JY,Z)-DJ(Y,JX,Z)=0\\
DJ(JZ,X,Y)-DJ(JX,Z,Y)+DJ(Z,JX,Y)-DJ(X,JZ,Y)=0\\
-DJ(JY,Z,X)+DJ(JZ,Y,X)-DJ(Y,JZ,X)+DJ(Z,JY,X)=0,
\end{\eq}
we obtain one direction and the other direction is obvious.
\qed
\begin{lem}\label{akc}
(almost K\"{a}hler condition) Let $(g,J,\omega)$ be an \aH structure. Then
\begin{\eq}
d\omega=0&\Longleftrightarrow& DJ(X,Y,Z)+DJ(Y,Z,X)+DJ(Z,X,Y)=0\\
&\Longrightarrow& DJ(JX,JY,Z)=-DJ(X,Y,Z)
\end{\eq}
\end{lem}
\pf
The first identity is directly obtained from definition since $D_{X}Y-D_{Y}X=[X,Y]$ and $DJ=D\omega$. For the second identity, adding the following identities together,
\begin{\eq}
DJ(JX,JY,JZ)+DJ(JY,JZ,JX)+DJ(JZ,JX,JY)&=&0\\
-DJ(JX,Y,Z)-DJ(Y,Z,JX)-DJ(Z,JX,Y)&=&0\\
DJ(X,JY,Z)+DJ(JY,Z,X)+DJ(Z,X,JY)&=&0\\
DJ(X,Y,JZ)+DJ(Y,JZ,X)+DJ(JZ,X,Y)&=&0
\end{\eq}
and by Lemma \ref{dj}, we obtain the desired result.
\qed

In \aH setting, modulo first order terms, $\triangle J$ is $(0,2)+(2,0)$.
\begin{lem}\label{lapj}
Let $(g,J,\omega)$ be an \aH structure. Then
\begin{\eq}
(\triangle J)^{(1,1)}=-\mathcal{N},~\text{where}\quad\mathcal{N}(X,Y)=DJ(i,j,JX)DJ(i,j,Y).
\end{\eq}
So $\triangle J+\mathcal{N}$ is $(0,2)+(2,0)$.
\end{lem}
\pf
\begin{eqnarray*}
\triangle J(JX,JY)&=&D_{i}(DJ(i,JX,JY))-DJ(i,j,JY)DJ(i,j,X)-DJ(i,JX,j)DJ(i,j,Y)\\
&=&-D_{i}(DJ(i,X,Y))+2DJ(i,j,JX)DJ(i,j,Y)\\
&=&-\triangle(X,Y)+2\mathcal{N}(X,Y).
\end{eqnarray*}

So we finish the proof.
\qed

\section{Deformation Conditions in Almost Hermitian Setting}\label{con}
First, we derive some algebraic conditions for the deformation of \aH structures, which is necessary.
\begin{lem}\label{01}
Let $(g,J)$ be a family of almost Hermitian structures, $\frac{\partial}{\partial t}g=h, \frac{\partial}{\partial t}J=K$. Then $(h,K)$ satisfies the following algebraic conditions.
\begin{\eq}
&(a)& h~ \text{is symmetric},\\
&(b)& K~ \text{is}~ (0,2)+(2,0),\\
&(c)& K^{sym}J=h^{(0,2)+(2,0)}.
\end{\eq}
\end{lem}
\pf
Condition (a) is obvious. For (b),
\begin{eqnarray*}
0=\frac{\partial}{\partial t}J^{2}=KJ+JK.
\end{eqnarray*}

For (c),
\begin{eqnarray*}
0&=&\frac{\partial}{\partial t}(g(JX,JY)-g(X,Y))\\
&=&h(JX,JY)-h(X,Y)+g(KX,JY)+g(JX,KY)\\
&=&-2h^{(0,2)+(2,0)}(X,Y)+K(JX,Y)+K(Y,JX)\\
&=&-2h^{(0,2)+(2,0)}(X,Y)+2(K^{sym}J)(X,Y).
\end{eqnarray*}
\qed
\begin{rem}\label{eta}
Similarly, we have
\begin{\eq}
\eta^{(0,2)+(2,0)}&=&K^{skew},\\
\eta^{(1,1)}&=&h^{(1,1)}J.
\end{\eq}
\end{rem}
\begin{lem}
Let $(g,J)$ be an almost Hermitian structure. Then $(L_{X}g,\:L_{X}J)$ satisfies the necessary condition of a deformation of $(g,J)$.
\end{lem}
\pf
Let $\phi_{t}$ be the 1-parameter transformation groups generated by $X$, $g_{t}=\phi_{t}^{\ast}g$, $J_{t}=\phi_{t}^{\ast}J$, then
\begin{eqnarray*}
\frac{\partial}{\partial t}\Big{|}_{t=0}g_{t}=L_{X}g,\quad
\frac{\partial}{\partial t}\Big{|}_{t=0}J_{t}=L_{X}J.
\end{eqnarray*}

Then the result follows from Lemma \ref{01}.
\qed

Next, we consider the analytic condition to deform \aH structures. We only consider the second order flows. And to ensure the short-time existence (on compact manifolds), we assume $(h,K)$ satisfies the following conditions,\\
(1) $(h,K)$ is of second order with respect to $(g,J)$,\\
(2) Modulo gauge transformation, the symbol of the linearization of $(h,K)$ is $|\xi|^{2}Id$, i.e., there exists a vector field $\overline{X}$, s.t,
\begin{\eq}
h+L_{\overline{X}}g&=&g^{ij}\partial_{i}\partial_{j}g+\mathcal{O}(\partial g,\partial J),\\
K+L_{\overline{X}}J&=&g^{ij}\partial_{i}\partial_{j}J+\mathcal{O}(\partial g,\partial J).
\end{\eq}
\begin{rem}
In the Ricci flow, $\overline{X}=g^{ab}(\Gamma_{ab}^{k}-\overline{\Gamma}_{ab}^{k})\frac{\partial}{\partial x^{k}}$, where $\overline{\Gamma}$ is the Christoffel symbol of a fixed background linear connection.
\end{rem}
\begin{rem}
If $(h,K)$ satisfies both the algebraic condition and the analytic condition, then we may apply the same method in \cite{ST symp} (see also \cite{Dai}) to show that there exist a family of \aH structures on a compact manifold for a short while, with given initial data.
\end{rem}
\begin{rem}
There are other second order ``canonical" curvature flows being constructed, but not satisfying our analytic condition, see \cite{LeWang}.
\end{rem}
\begin{rem}
Naturally, we also require that our flow should coincide with the K\"{a}hler Ricci flow if the initial data is K\"{a}hler. But in fact, in K\"{a}hler setting, the K\"{a}hler Ricci flow is the unique flow satisfies the conditions above. So this requirement is vacant.
\end{rem}
From the above discussion, we see to classify the second order curvature flows, we just need to do the following two steps.
\\
Step 1: Classify the tensors up to order $2$, as well as gauge terms. More precisely, we need to classify the tensors of the following type,
\begin{\eq}
\text{first order terms: } \partial J*\partial J,\quad \text{second order terms: } \partial^{2} J,~ \partial^{2} g,\quad \text{gauge terms: }\overline{X}.
\end{\eq}
Step 2: Calculate the symbols of the second order terms and the gauge terms, and then find the suitable tensors to satisfy the analytic condition.
\begin{rem}
We only consider the quadratic terms of $\partial J$, for the consideration of scaling property.
\end{rem}


\section{Classification of Lower Order Tensors}\label{tensor}
We only consider the natural tensors, i.e., defined from $(g,J,\omega)$. For the classification, we mean that we give a list of tensors which satisfy some conditions, and all the tensors satisfying this condition are the linear combination of the tensors in the list. For a $2$-tensor $T$, ``modulo the transposition action" means that we regard $T$ and ${}^{t}T$ as the same tensor. And ``modulo the $J$ action" means that we regard $T$ and $J*T$ as the same tensor.

We only consider Levi-Civita connection since the difference between two connections also gives a tensor.

First, we consider the first order tensors. We require $h$ is of even type and $K$ is of odd type. So we only consider tensors of even or odd type. For instance, 
we don't consider tensors like $g+\omega$. So the zero order tensors are only $g$ and $\omega$. And to take contraction, we can only use $g$ and $\omega$.

We notice that in \aH setting, last two variables of $DJ$ is $(0,2)+(2,0)$. So it will vanish if we either take trace or complex trace in last two positions.
\begin{lem}\label{djdj}
Consider $2$-tensors of even type of form $DJ*DJ$. We take twice trace or complex trace of $DJ(\cdot,\cdot,\cdot)DJ(\cdot,\cdot,\cdot)$, then there are two positions remaining for the variables. Modulo the transposition and $J$ action, the tensors described above can be classified as follows,
\begin{\eq}
&&B^{1}(X,Y)=DJ(X,i,j)DJ(Y,i,j),\\
&&B^{2}(X,Y)=DJ(i,X,j)DJ(i,Y,j),\\
&&B^{3}(X,Y)=DJ(i,X,j)DJ(j,Y,i),\\
&&B^{4}(X,Y)=DJ(X,i,j)DJ(i,Y,j),\\
&&B^{5}(X,Y)=DJ(i,X,i)DJ(j,Y,j),\\
&&B^{6}(X,Y)=DJ(X,Y,i)DJ(j,i,j),\\
&&B^{7}(X,Y)=DJ(i,X,Y)DJ(j,i,j),\\
&&B^{8}(X,Y)=DJ(JX,i,j)DJ(Y,Ji,j),\\
&&B^{9}(X,Y)=DJ(i,JX,j)DJ(Ji,Y,j),\\
&&B^{10}(X,Y)=DJ(i,JX,Y)DJ(j,Ji,j).
\end{\eq}
Furthermore $B^1, B^3, B^5$ are symmetric, $B^2, B^{9}$ are symmetric and $(1,1)$, $B^7, B^{10}$ are skew and $(0,2)+(2,0)$, and ${}^{t}B^{8}=J^{*}B^{8}$.
\end{lem}
\pf
We identify $(X,Y)$, $(Y,X)$, $(JY,JX)$ and $(JY,JX)$. And we frequently use Lemma \ref{dj} implicitly. First, we consider the case that the variables are $X$, $Y$, then we need to take twice trace or twice complex trace in remaining four positions. Suppose $X$ is in the first position, i.e.
\begin{\eq}
DJ(X,\cdot,\cdot)DJ(\cdot,\cdot,\cdot).
\end{\eq}

There are three ways to pose $Y$,
\begin{\eq}
DJ(X,\cdot,\cdot)DJ(Y,\cdot,\cdot), \quad DJ(X,\cdot,\cdot)DJ(\cdot,Y,\cdot), \quad DJ(X,Y,\cdot)DJ(\cdot,\cdot,\cdot).
\end{\eq}

We obtain $B^{1}$, $B^{4}$, $B^{6}$ respectively. If $X$ (and $Y$) is not in the first position, we may assume $X$ is in the second position. There are two positions for $Y$,
\begin{\eq}
DJ(\cdot,X,\cdot)DJ(\cdot,Y,\cdot), \quad DJ(\cdot,X,Y)DJ(\cdot,\cdot,\cdot).
\end{\eq}

If we take twice trace of $DJ(\cdot,X,\cdot)DJ(\cdot,Y,\cdot)$, it gives $B^{2}$, $B^{3}$ and $B^{5}$. And if we take twice complex trace of $DJ(\cdot,X,\cdot)DJ(\cdot,Y,\cdot)$, it gives $B^{9}$. We obtain $B^{7}$ from $DJ(\cdot,X,Y)DJ(\cdot,\cdot,\cdot)$. Next, we consider the case that the variables are $JX$ and $Y$. Then we need to take once trace and once complex trace in the remaining four positions. If $JX$ is in the first position, we only need to consider
\begin{\eq}
DJ(JX,\cdot,\cdot)DJ(Y,\cdot,\cdot),
\end{\eq}

since other cases are reduced to the above situation. Then we obtain $B^{8}$. If $JX$ is in the second position, we only need to consider
\begin{\eq}
DJ(\cdot,JX,Y)DJ(\cdot,\cdot,\cdot),
\end{\eq}
for the same reason. Then we obtain $B^{10}$. From the definition, we can easily obtain the identities of the transposition and $J$ action of $B^{1}$ to $B^{10}$.
\qed
\begin{lem}
The functions of form $DJ*DJ$ can be classified as follows,
\begin{\eq}
E^{1}&=&DJ(i,j,k)DJ(i,j,k)=|DJ|^{2},\\
E^{2}&=&DJ(i,j,k)DJ(j,i,k),\\
E^{3}&=&DJ(i,i,j)DJ(k,k,j)=|\theta|^{2},\\
E^{4}&=&DJ(i,j,k)DJ(Ji,Jj,k).
\end{\eq}
\end{lem}
\pf
We just take trace or complex trace of $B^{1}$ to $B^{10}$ to obtain the desired functions. Notice that for a symmetric $2$-tensor, it will vanish if we take complex trace. And for $B^{7}$ and $B^{10}$, they vanish either we take trace or complex trace. For $B^{4}$,
\begin{\eq}
DJ(k,i,j)DJ(i,Jk,j)=DJ(k,i,j)DJ(i,k,Jj)=-DJ(i,k,Jj)DJ(k,i,j).
\end{\eq}

So the complex trace of $B^{4}$ vanishes. The similar reason hold for $B^{6}$, $B^{8}$, $B^{9}$. Then we obtain $E^{1}$ (from $B^{1},B^{2}$), $E^{2}$ (from $B^{3},B^{4}$), $E^{4}$ (from $B^{5},B^{6}$), $E^{3}$ (from $B^{8},B^{9}$).
\qed

To summarize
\begin{lem}
Modulo the transposition and $J$ action, the $2$-tensors of form $DJ*DJ$ can be classified as follows,
\begin{\eq}
\text{Even type}:&& \quad B^{i},~ 1\leq i \leq 10,\qquad E^{i}g, ~ 1\leq i \leq 4,\\
\text{Odd type}:&& \quad B^{i}J,~ 1\leq i \leq 10,\qquad E^{i}\omega, ~ 1\leq i \leq 4.
\end{\eq}
Functions of form $DJ*DJ$ can be classified as $E^{i}$, $1\leq i \leq 4.$
\end{lem}

Now we consider the second order tensors. To begin with, we consider the tensors in terms of $\partial^{2}g$, but without $\partial^{2}J$. They are tensors in terms of the Riemannian curvature, more precisely, they are of the form $Rm*\underbrace {J*\dots*J}_{m \text{th}}$, for any $m\geq 0$.
\begin{lem}\label{cur}
Modulo the transposition and $J$ action, the $2$-tensors in terms of the Riemannian curvature can be classified as follows,
\begin{\eq}
\text{Even type:}&&\quad Ric(X,Y),\quad \rho'(JX,Y),\quad Rg,\quad s'g.\\
\text{Odd type:}&&\quad Ric(JX,Y),\quad \rho'(X,Y),\quad R\omega,\quad s'\omega.
\end{\eq}
The functions in terms of the Riemmanian curvature are given by $R$ and $s'$.
\end{lem}
\pf
We only consider the case that the tensor is of even type, and $Rm$ is not totally contracted. From the symmetries of the Riemmanian curvature, we see there are three possibilities, modulo the transposition and $J$ action,
\begin{\eq}
Ric(X,Y)=Rm(X,i,i,Y),\quad Rm(JX,i,Ji,Y),\quad\rho'(X,Y)=Rm(JX,Y,i,Ji).
\end{\eq}

But from Bianchi identity,
\begin{\eq}
Rm(JX,i,Ji,Y)&=&-Rm(JX,Ji,Y,i)-Rm(JX,Y,i,Ji)\\
&=&-Rm(JX,i,Ji,Y)-Rm(JX,Y,i,Ji).
\end{\eq}

So
\begin{eqnarray*}
Rm(JX,i,Ji,Y)=-\frac{1}{2}\rho'(X,Y).\label{rp}
\end{eqnarray*}

Then we can obtain the result.
\qed

Now we consider the tensors in terms of $\partial^{2}J$, more precisely, they are of the form $D^2J*\underbrace {J*\dots*J}_{m \text{th}}$, for any $m\geq 0$.. Notice that if we take trace or complex trace in two positions of $D^{2}J$, it will reduce to the first order terms. For instance,
\begin{\eq}
D^{2}J(X,Y,i,Ji)=D_{X}(DJ(Y,i,Ji))-DJ(Y,i,D_{X}(Ji))=-DJ(Y,i,j)DJ(X,i,j).
\end{\eq}
\begin{lem}\label{d2j}
Modulo the transposition action and $J$ action, modulo the lower order terms and Riemannian curvature terms, the $2$-tensors in terms of $\partial^{2}J$ can be classified as follows
\begin{\eq}
\text{Odd type:}&&\triangle J(X,Y),\quad D^{2}J(X,i,Y,i),\\
\text{Even type:}&&\triangle J(JX,Y),\quad D^{2}J(JX,i,Y,i).
\end{\eq}
There is no functions in terms of $\partial^{2}J$.
\end{lem}
\pf
Notice that modulo the Riemannian curvature, the first two components of $D^{2}J(\cdot,\cdot,\cdot,\cdot)$ is symmetric. And modulo the lower order terms in $DJ$, the last two components of $D^{2}J(\cdot,\cdot,\cdot,\cdot)$ is skew and $(0,2)+(2,0)$. So suppose the tensors are of odd type, and $D^{2}J$ is not totally contracted, then the tensors are classified as
\begin{\eq}
\triangle J(X,Y),\quad D^{2}J(X,i,Y,i).
\end{\eq}
Now we consider the contraction of above tensors. From Lemma \ref{lapj}, we see modulo lower order tensors, $\triangle J$ is $(0,2)+(2,0)$. And notice that $\triangle J=\triangle \omega$ is skew. So either taking trace or complex trace of $\triangle J$ will vanish. For $D^{2}J(j,i,j,i)$ and $D^{2}J(j,i,Jj,i)$, notice that, modulo the Riemannian curvature and modulo the lower order terms, the $i,j$ in the first two positions commutes and $i,j$ in the last two positions anti-commutes. So they also vanish. Hence there is no functions in terms of $\partial^{2}J$.
\qed

We finish the classification of the lower order tensors in \aH setting. Now we consider two special cases. We only state the results which can be reduced from the results in the setting of \aH. First, we consider \aK setting.

We notice that in \aK setting, every two variables of $DJ$ is $(0,2)+(2,0)$. So it will vanish if we either take trace or complex trace in any two positions. And if we take trace or complex trace in any two of the last three positions of $D^{2}J$, it will reduce to the first order terms.
\begin{lem}\label{akt}
In \aK setting, the $2$-tensors of form $DJ*DJ$ can be classified as
\begin{\eq}
B^{1},\quad B^{2},\quad|DJ|^2g,\quad B^{1}J,\quad B^{2}J,\quad|DJ|^2\omega.
\end{\eq}

Furthermore, both $B^{1}$ and $B^{2}$ are symmetric and $(1,1)$.

The functions of form $DJ*DJ$ can be classified as $|DJ|^{2}$.
\end{lem}
\pf
We just need to show the list in Lemma \ref{djdj} can be reduced to $B^{1}$, $B^{2}$. From Lemma \ref{akc}, we see $B^{5},B^{6},B^{7},B^{10}$ vanish, and $B^{8},B^{9}$ are reduced to $B^{1},B^{2}$. So we only need express $B^{3},B^{4}$ in terms of $B^{1},B^{2}$. In fact,
\begin{\eq}
DJ(X,i,j)DJ(i,Y,j)&=&-DJ(X,i,j)(DJ(Y,j,i)+DJ(j,i,Y))\\
&=&DJ(X,i,j)DJ(Y,i,j)-DJ(X,j,i)DJ(j,Y,i).
\end{\eq}

So $B^{4}=\frac{1}{2}B^{1}$. And
\begin{\eq}
DJ(i,X,j)DJ(j,Y,)&=&-DJ(i,X,j)(DJ(Y,i,j)+DJ(i,j,Y))\\
&=&-DJ(Y,i,j)DJ(i,X,j)+DJ(i,X,j)DJ(i,Y,j).
\end{\eq}

So $B^{3}=-\frac{1}{2}B^{1}+B^{2}$. And it is to see, in the \aK setting, both $B^{1}$ and $B^{2}$ are symmetric and $(1,1)$. So we finish the proof.
\qed
\begin{lem}
In \aK setting, $2$-tensors in terms of $\partial^{2}J$ can be classifies as follows
\begin{\eq}
\triangle J(X,Y), \quad \triangle J(JX,Y).
\end{\eq}
\end{lem}
\pf
Notice that $D^{2}J(X,i,Y,i)$ is reduced to the first order terms. (In fact, it vanishes.)
\qed
\begin{lem}
In \aK setting, the second order functions are classified as $R$, $|DJ|^{2}$.
\end{lem}
The following lemma gives the proof.
\begin{lem}
In \aK setting, We have the following identity.
\begin{\eq}
s'+2R+|DJ|^{2}=0.
\end{\eq}
\end{lem}
\pf
First we have
\begin{eqnarray*}
Rm(i,j,Jk,l)+Rm(i,j,k,Jl)=D^{2}J(i,j,k,l)-D^{2}J(j,i,k,l).\label{ric}
\end{eqnarray*}

To see this, suppose $i,j,k,l$ are in a local normal coordinate chart. Then
\begin{\eq}
D_{i}D_{j}(Jk)=D_{i}(DJ(j,k)+JD_{j}k)=D^{2}J(i,j,k)+JD_{i}D_{j}k
\end{\eq}

So
\begin{\eq}
Rm(i,j,Jk,l)&=&\langle D_{i}D_{j}(Jk)-D_{j}D_{i}(Jk),l\rangle\\
&=&D^{2}J(i,j,k,l)-D^{2}J(j,i,k,l)-\langle D_{i}D_{j}k,Jl\rangle+\langle D_{j}D_{i}k,Jl\rangle\\
&=&D^{2}J(i,j,k,l)-D^{2}J(j,i,k,l)-Rm(i,j,k,Jl)
\end{\eq}

So we obtain the identity. (One may also show it directly from Ricci identity.)\\
Now let $k=j$, $l=Ji$, we see
\begin{\eq}
Rm(i,j,Jj,Ji)-R=D^{2}J(i,j,j,Ji)-D^{2}J(j,i,j,Ji).
\end{\eq}

From $d\omega=0$,
\begin{\eq}
-D^{2}J(j,i,j,Ji)&=&-D_{j}(DJ(i,j,Ji))+DJ(i,j,k)DJ(j,i,k)\\
&=&D_{j}(DJ(j,Ji,i)+DJ(Ji,i,j))+trB^{3}\\
&=&D_{j}(-DJ(j,i,Ji)-DJ(i,i,Jj))+\frac{1}{2}|DJ|^{2}\\
&=&-D^{2}J(j,j,i,Ji)-|DJ|^{2}-D^{2}J(j,i,i,Jj)-|\theta|^{2}+\frac{1}{2}|DJ|^{2}.
\end{\eq}

Since $DJ(i,i,X)=0$, we have
\begin{\eq}
Rm(i,j,Jj,Ji)-R=-\triangle J(i,Ji)-\frac{1}{2}|DJ|^{2}.
\end{\eq}

From Lemma \ref{lapj}, we see
\begin{\eq}
\triangle J(i,Ji)=-|DJ|^{2}.
\end{\eq}

And from the identity in (\ref{cur}), we see
\begin{\eq}
Rm(i,j,Jj,Ji)=-\frac{1}{2}s'
\end{\eq}

So we finish the proof.
\qed

In the case of dimension $4$, we can reduce the tensors further in \aK setting. The author learned the well-known result below from R. Bryant in mathoverflow\cite{Bryant}.
\begin{lem}
Let $(M,g,J,\omega)$ be an \aK manifold with $\dim_{\mathbb{R}}M=4$. Let $p\in M$. Then in any local coordinate chart of $p$, there exists a local unitary frame, i.e.,
\begin{displaymath}
g=\left(
\begin{array}{cccc}
1 & & &\\
& 1 & &\\
& & 1 &\\
& & & 1
\end{array}\right),\quad
J=\left(
\begin{array}{cccc}
& 1 & &\\
-1 & & &\\
& & & 1\\
& & -1 &
\end{array}\right),
\end{displaymath}
such that at $p$,
\begin{displaymath}
B^1=\left(
\begin{array}{cccc}
4a^2 & & &\\
& 4a^2 & &\\
& & 0 &\\
& & & 0
\end{array}\right),
\quad
B^2=\left(
\begin{array}{cccc}
2a^2 & & &\\
& 2a^2 & &\\
& & 2a^2 &\\
& & & 2a^2
\end{array}\right),
\end{displaymath}
where $a^2=\frac{1}{8}|DJ|^2$. In particular, $B^{2}=\frac{1}{4}|DJ|^{2}g$.
\end{lem}
\pf
We can choose a local coframe $\{\eta_{i}\},i=1,2,3,4$, s.t.
\begin{\eq}
g=\eta_{1}^{2}+\eta_{2}^{2}+\eta_{3}^{2}+\eta_{4}^{2},\quad \omega=\eta_{1}\wedge\eta_{2}+\eta_{3}\wedge\eta_{4}.
\end{\eq}

Notice that $\omega$ is self-dual. Let
\begin{\eq}
\kappa=\eta_{1}\wedge\eta_{3}+\eta_{4}\wedge\eta_{2},\quad \lambda=\eta_{1}\wedge\eta_{4}+\eta_{2}\wedge\eta_{3}.
\end{\eq}

Then $(\omega,\kappa,\lambda)$ forms a pointwise basis for the self-dual $2$-forms. First we claim $D_{X}\omega$ is self-dual. In fact, in general, we have
\begin{\eq}
D_{X}(*\alpha)=*D_{X}\alpha.
\end{\eq}

This is from
\begin{\eq}
D_{X}(*\alpha)\wedge \beta&=&D_{X}(*\alpha\wedge \beta)-*\alpha\wedge D_{X}\beta=D_{X}(\langle \alpha,\beta\rangle dV)-\langle \alpha,D_{X}\beta\rangle dV\\
&=&\langle D_{X}\alpha,\beta\rangle dV=*D_{X}\alpha\wedge \beta.
\end{\eq}

So set
\begin{\eq}
D\omega=\alpha\otimes\kappa+\beta\otimes\lambda+\gamma\otimes\omega,
\end{\eq}

where $\alpha,\beta,\gamma$ are $1$-forms.
But notice that
\begin{\eq}
D_{X}\omega\wedge\omega=D_{X}(dV)=0,
\end{\eq}

so $\gamma=0$. Hence we have
\begin{\eq}
0=d\omega=\alpha\wedge\kappa+\beta\wedge\lambda.
\end{\eq}

Let
\begin{\eq}
\alpha=a_{1}\eta_{1}+a_{2}\eta_{2}+a_{3}\eta_{3}+a_{4}\eta_{4}.
\end{\eq}

Then
\begin{\eq}
\beta=a_{2}\eta_{1}-a_{1}\eta_{2}+a_{4}\eta_{3}-a_{3}\eta_{4}.
\end{\eq}

If $\alpha(p)=0$, then there is nothing to prove. If $\alpha(p)\neq 0$, then we may choose
\begin{\eq}
a_{1}(p)=a>0, ~a_{2}(p)=a_{3}(p)=a_{4}(p)=0.
\end{\eq}

Hence at $p$,
\begin{\eq}
DJ=a\eta_{1}\otimes(\eta_{1}\wedge\eta_{3}+\eta_{4}\wedge\eta_{2})-a\eta_{2}\otimes(\eta_{1}\wedge\eta_{4}+\eta_{2}\wedge\eta_{3}).
\end{\eq}

So at $p$,
\begin{\eq}
D_{1}J=a\kappa,~D_{2}J=-a\lambda,~D_{3}J=D_{4}J=0.
\end{\eq}

Therefore at $p$,
\begin{displaymath}
B^1=\left(
\begin{array}{cccc}
4a^2 & & &\\
& 4a^2 & &\\
& & 0 &\\
& & & 0
\end{array}\right),
\quad
B^2=\left(
\begin{array}{cccc}
2a^2 & & &\\
& 2a^2 & &\\
& & 2a^2 &\\
& & & 2a^2
\end{array}\right).
\end{displaymath}

Taking trace, we see $a^{2}=\frac{1}{8}|DJ|^{2}$.
\qed
\begin{rem}
If $DJ\neq 0$ at $p$, we can choose a local unitary frame such that the above result holds in a neighborhood of $p$, not only at $p$.
\end{rem}
\begin{rem}
From the proposition above, we see ~$\frac{1}{2}B^{1}-B^{2}\leq 0$. So the Symplectic Curvature Flow is ``slower" than Ricci Flow in $g$.
\end{rem}

Now we consider Hermitian setting.
\begin{lem}\label{tenh}
In Hermitian setting, the $2$-tensors of form $DJ*DJ$ can be classified as
\begin{\eq}
B^{i},1\leq i\leq3,5\leq i\leq7, ~E^{i}g,1\leq i\leq 3,\quad B^{i}J,1\leq i\leq3,5\leq i\leq7, ~E^{i}\omega,1\leq i\leq 3.
\end{\eq}
Furthermore $B^{1},B^{2}$ are symmetric and $(1,1)$, $B^3$ is symmetric and $(0,2)+(2,0)$, $B^{5}$ is symmetric, $B^{6}$ is $(1,1)$ and $B^7$ is skew and $(0,2)+(2,0)$.\\
The functions of type $DJ*DJ$ can be classified as $E^{i},1\leq i\leq 3$.
\end{lem}
\pf
By taking advantage of Lemma \ref{hc}, the tensors in the list of Lemma \ref{djdj} are reduced to the result above.
\qed

Next, we consider the gauge terms. More precisely, we consider the vector fields defined from $(g,J,\omega,\overline{D})$, where $\overline{D}$ is a fixed background linear connection. Naturally, we require the vector field is of first order and of even type. First, we have ``canonical" gauge.
\begin{lem}
Let $X$ be a first order vector field defined from $(g,J,\omega)$, then $X$ can be classified as
\begin{\eq}
\theta^{\sharp}=DJ(i,Ji),\quad J\theta^{\sharp}.
\end{\eq}
\end{lem}
\pf
We just need take once trace or once complex trace of $DJ$. It will vanish if we contract the last two components, so we have the result above.
\qed

Now, we consider the vector field in terms of $\overline{D}$, i.e., $\overline{D}g$ and $\overline{D}J$. Notice that
\begin{\eq}
\overline{D}J-DJ=(\overline{D}-D)*J.
\end{\eq}

So we only need to classify the vector field in terms of $\overline{D}g$. Notice that the last two components of $\overline{D}g$ is symmetric, so we have the following result.
\begin{lem}
The gauge terms in $(\partial g, J)$ can be classified as
\begin{\eq}
\overline{X_{1}}^{\flat}=g^{ij}\overline{D}g(i,j,\cdot), \quad\overline{X_{2}}^{\flat}=g^{ij}\overline{D}g(\cdot,i,j), \quad\overline{X_{0}}^{\flat}=g^{ij}J_{j}^{k}\overline{D}g(i,k,J\cdot).
\end{\eq}
\end{lem}
\begin{rem}
In Ricci Flow, the gauge term is $\overline{X}=\overline{X_{1}}-\frac{1}{2}\overline{X_{2}}.$
\end{rem}
We can calculate $L_{\theta^{\sharp}}J$ as follows.
\begin{lem}\label{lee}
Let $(g,J)$ be an almost Hermitian structure, then
\begin{\eq}
L_{Z}J(X,Y)&=&DJ(Z,X,Y)-DZ(JX,Y)-DZ(X,JY).\\
L_{\theta^{\sharp}}J(X,Y)&=&-D^{2}J(JX,i,Ji,Y)-D^{2}J(X,i,Ji,JY)\\
&&-DJ(JX,i,j)DJ(i,j,Y)-DJ(X,i,j)DJ(i,j,JY)+DJ(j,X,Y)DJ(i,Ji,j).
\end{\eq}
\end{lem}
\pf
For the first one,
\begin{\eq}
L_{Z}J(X)&=&[Z,JX]-J[Z,X]=D_{Z}(JX)-D_{JX}Z-JD_{Z}X+JD_{X}Z\\
&=&DJ(Z,X)-DZ(JX)+JDZ(X).
\end{\eq}

So we obtain the first identity. For the second one,
\begin{\eq}
D\theta^{\sharp}(X,Y)=g(D_{X}(DJ(i,Ji)),Y)=D^{2}J(X,i,Ji,Y)+DJ(i,j,Y)DJ(X,i,j),
\end{\eq}

then we finish the proof.
\qed

\section{Classification of Second Order Curvature Flows}\label{flow}
To classify the second order curvature flows, we just need to find the suitable tensor satisfying the conditions in Section \ref{con}. We have discussed the algebraic condition in Section \ref{tensor}. Now we consider the analytic condition which means we need to calculate the symbols of the second order tensors. The result of symbol calculations in non-K\"{a}hler geometry is well-known in literature. For the convenience of readers, we still do calculations here.\\
\textit{Proof of Theorem \ref{ah}:}
First, we consider the evolution of $J$, $\delta J=K$. From our classification result Lemma \ref{cur} and Lemma \ref{d2j}, and considering the parity of  the type of $J$, we see the candidates of the second order terms are
\begin{\eq}
&&\triangle J(X,Y),~ D^{2}J(X,i,Y,i),~ D^{2}J(JX,i,JY,i),~ D^{2}J(Y,i,X,i),~ D^{2}J(JY,i,JX,i),\\
&&Ric(JX,Y),~ Ric(X,JY),~ \rho'(JX,Y),~ \rho'(X,JY),~ RJ,~ s'J.
\end{\eq}

Here $\triangle J$ and $\rho'(J\cdot,\cdot)$ are skew and $\mathcal{R}$ is symmetric, so we don't need to consider their transposition.
Considering algebraic condition, the tensors of evolution terms should be $(0,2)+(2,0)$ and from Lemma \ref{lee} we see $D^{2}J(X,i,Y,i)+D^{2}J(JX,i,JY,i)$ is just $L_{\theta^{\sharp}}J$ modulo lower order terms. So our candidates of the second order terms are in fact
\begin{\eq}
\triangle J,~{}^{t}(L_{\theta^{\sharp}}J),~\mathcal{R},~ (\rho'(J\cdot,\cdot))^{(0,2)+(2,0)}.
\end{\eq}

Now we calculate their symbols.
\begin{lem}
In local coordinate chart,
\begin{\eq}
Rm_{ijkl}&=&\frac{1}{2}(\partial_{i}\partial_{k}g_{lj}-\partial_{i}\partial_{l}g_{jk}
-\partial_{j}\partial_{k}g_{li}+\partial_{j}\partial_{l}g_{ik})+\mathcal{O}(\partial g).\\
(D^{2}J)_{ijk}^{l}&=&\partial_{i}\partial_{j}J_{k}^{l}
+\frac{1}{2}J_{k}^{p}g^{lq}(\partial_{i}\partial_{j}g_{qp}+\partial_{i}\partial_{p}g_{qj}-\partial_{i}\partial_{q}g_{jp})\\
&&-\frac{1}{2}J_{p}^{l}g^{pq}(\partial_{i}\partial_{j}g_{qk}+\partial_{i}\partial_{k}g_{qj}-\partial_{i}\partial_{q}g_{jk})+\mathcal{O}(\partial g,\partial J).
\end{\eq}
\end{lem}
\pf
In any local coordinate chart,
\begin{\eq}
Rm_{ijkl}&=&\partial_{i}\Gamma_{jkl}-\partial_{j}\Gamma_{ikl}+\mathcal{O}(\partial g)\\
&=&\frac{1}{2}(\partial_{i}\partial_{j}g_{lk}+\partial_{i}\partial_{k}g_{lj}-\partial_{i}\partial_{l}g_{jk})
-\frac{1}{2}(\partial_{j}\partial_{i}g_{lk}+\partial_{j}\partial_{k}g_{li}-\partial_{j}\partial_{l}g_{ik})+\mathcal{O}(\partial g)\\
&=&\frac{1}{2}(\partial_{i}\partial_{k}g_{lj}-\partial_{i}\partial_{l}g_{jk}-\partial_{j}\partial_{k}g_{li}+\partial_{j}\partial_{l}g_{ik})+\mathcal{O}(\partial g).\\
(D^{2}J)_{ijk}^{l}&=&(D_{i}(D_{j}J(k)))^{l}+\mathcal{O}(\partial g,\partial J)\\
&=&\partial_{i}\partial_{j}J_{k}^{l}+J_{k}^{p}\partial_{i}\Gamma_{jp}^{l}-J_{p}^{l}\partial_{i}\Gamma_{jk}^{p}+\mathcal{O}(\partial g,\partial J)\\
&=&\partial_{i}\partial_{j}J_{k}^{l}+\frac{1}{2}J_{k}^{p}g^{lq}(\partial_{i}\partial_{j}g_{qp}+\partial_{i}\partial_{p}g_{qj}-\partial_{i}\partial_{q}g_{jp})\\
&&-\frac{1}{2}J_{p}^{l}g^{pq}(\partial_{i}\partial_{j}g_{qk}+\partial_{i}\partial_{k}g_{qj}-\partial_{i}\partial_{q}g_{jk})+\mathcal{O}(\partial g,\partial J).
\end{\eq}
\qed

To calculate the symbols of our candidates, we notice that symbol is also tensorial, so we just need to calculate the symbols of $Rm$ and $D^{2}J$, and then use the corresponding manner to take trace to get the desired symbols. To calculate the symbols of the linearization operators of $Rm$ and $D^{2}J$, what we need to do is to replace the second derivative terms by their deformation terms, and replace ``$\partial_{i}$" by $\xi_{i}$, where $\xi$ is a $1$-form. Finally, we simplified the tensors we obtained. We denote $\sigma$ to be the symbol of the linearization operator of a tensor. Then for instance
\begin{\eq}
\sigma: g^{ij}\partial_{i}\partial_{j}g_{ab}\mapsto g^{ij}\xi_{i}\xi_{j}h_{ab}=|\xi|^{2}h_{ab}.
\end{\eq}

Since symbol is also tensorial, we also use the Riemannian metric to identify $TM$ and $T^{*}M$. We may use orthonormal frame to reduce the calculation, and notice that $J_{i}^{j}=-J_{j}^{i}$. Then we obtain the following results.
\begin{lem}\label{orisym}
\begin{\eq}
\sigma(Rm)_{ijkl}&=&\frac{1}{2}(\xi_{i}\xi_{k}h_{lj}-\xi_{i}\xi_{l}h_{jk}-\xi_{j}\xi_{k}h_{li}+\xi_{j}\xi_{l}h_{ik}),\\
\sigma(D^{2}J)_{ijkl}&=&\xi_{i}\xi_{j}K_{kl}+\frac{1}{2}(\xi_{i}\xi_{j}J_{k}^{p}h_{lp}+\xi_{i}\xi_{p}J_{k}^{p}h_{lj}-\xi_{i}\xi_{l}J_{k}^{p}h_{jp}\\
&&+\xi_{i}\xi_{j}J_{l}^{p}h_{pk}+\xi_{i}\xi_{k}J_{l}^{p}h_{pj}-\xi_{i}\xi_{p}J_{l}^{p}h_{jk}).
\end{\eq}
\end{lem}
\begin{lem}\label{symbol}
\begin{\eq}
\sigma(\triangle J)_{ab}&=&|\xi|^{2}K(a,b)+\frac{1}{2}(|\xi|^{2}h(Ja,b)+h(\xi,b)\xi(Ja)-h(\xi,Ja)\xi(b)\\
&&+|\xi|^{2}h(Jb,a)+h(\xi,Jb)\xi(a)-h(\xi,a)\xi(Jb)),\\
\sigma({}^{t}(L_{\theta^{\sharp}}J))_{ab}&=&K(a,\xi)\xi(b)+\frac{1}{2}\xi(b)\xi(Ja)\text{tr}h+h(J\xi,a)\xi(b)\\
&&-K(Ja,\xi)\xi(Jb)+\frac{1}{2}\xi(Jb)\xi(a)\text{tr}h-h(J\xi,Ja)\xi(Jb),\\
\sigma (\mathcal{R})_{ab}&=&\frac{1}{2}(h(b,\xi)\xi(Ja)+h(Ja,\xi)\xi(b)-|\xi|^{2}h(Ja,b)-\xi(Ja)\xi(b)trh\\
&&+h(Jb,\xi)\xi(a)+h(a,\xi)\xi(Jb)-|\xi|^{2}h(a,Jb)-\xi(a)\xi(Jb)trh),\\
\sigma((\rho'(J\cdot,\cdot))^{(0,2)+(2,0)})_{ab}&=&h(a,J\xi)\xi(b)-h(b,J\xi)\xi(a)-h(Ja,J\xi)\xi(Jb)+h(Jb,J\xi)\xi(Ja).
\end{\eq}
\end{lem}
For the gauge terms, the candidates are
\begin{\eq}
L_{\theta^{\sharp}}J,~L_{\overline{X_{1}}}J,~L_{\overline{X_{2}}}J,~L_{\overline{X_{3}}}J.
\end{\eq}
We have calculated $\sigma({}^{t}(L_{\theta^{\sharp}}J))$. For other symbols, by definition, we have
\begin{lem}
\begin{\eq}
\overline{X_{1}}^{k}&=&g^{kl}g^{ij}\partial_{i}g_{jl}+\mathcal{O}(g),\\
\overline{X_{2}}^{k}&=&g^{kl}g^{ij}\partial_{l}g_{ij}+\mathcal{O}(g),\\
\overline{X_{0}}^{k}&=&g^{jk}J_{j}^{l}g^{ab}J_{b}^{i}\partial_{a}g_{il}+\mathcal{O}(g, J).\\
\end{\eq}
\end{lem}
From the Lemma \ref{lee}, we have
\begin{lem}
\begin{\eq}
(L_{\overline{X_{1}}}J)_{a}^{b}&=&-J_{a}^{k}g^{bl}g^{pq}\partial_{k}\partial_{p}g_{ql}+J_{k}^{b}g^{kl}g^{pq}\partial_{a}\partial_{p}g_{ql}
+\mathcal{O}(\partial g,\partial J),\\
(L_{\overline{X_{2}}}J)_{a}^{b}&=&-J_{a}^{k}g^{bl}g^{pq}\partial_{k}\partial_{l}g_{pq}+J_{k}^{b}g^{kl}g^{pq}\partial_{a}\partial_{l}g_{pq}
+\mathcal{O}(\partial g,\partial J),\\
(L_{\overline{X_{0}}}J)_{a}^{b}&=&-J_{a}^{k}g^{bp}J_{p}^{l}g^{ij}J_{j}^{q}\partial_{k}\partial_{i}g_{ql}+g^{bl}g^{ij}J_{j}^{q}\partial_{a}\partial_{i}g_{ql}
+\mathcal{O}(\partial g,\partial J).
\end{\eq}
\end{lem}
Then the symbols of the above tensors are
\begin{lem}
\begin{\eq}
\sigma(L_{\overline{X_{1}}}J)_{ab}=-\xi(Ja)h(\xi,b)-\xi(a)h(\xi,Jb),\\
\sigma(L_{\overline{X_{2}}}J)_{ab}=-\xi(Ja)\xi(b)trh-\xi(a)\xi(Jb)trh,\\
\sigma(L_{\overline{X_{0}}}J)_{ab}=-\xi(Ja)h(J\xi,Jb)+\xi(a)h(J\xi,b).
\end{\eq}
\end{lem}
We consider the symbol of $K$ with respect to $J$ first. We assume there is no $L_{\theta^{\sharp}}J$ in $K$. Then we only need to consider $\triangle J$ and ${}^{t}L_{\theta^{\sharp}}J$. Notice that the symbol with respect to $J$ of $\triangle J$ is already good and we cannot compensate the symbol with respect to $J$ in ${}^{t}L_{\theta^{\sharp}}J$ by using $L_{\theta^{\sharp}}J$ (other gauge terms only involving $\partial g$), so we can only choose $\triangle J$ to be the symbol term of $K$. As for the symbol of $K$ with respect to $h$, we see that to compensate terms of $|\xi|^{2}hJ$ in the symbol of $\triangle J$, we must have $\triangle J+\mathcal{R}$, and notice that $\triangle J+\mathcal{R}+L_{\overline{X}}J$ gives us desirable symbol. As for $(\rho'(J\cdot,\cdot))^{(0,2)+(2,0)}$, we can't compensate $h(Ja,J\xi)\xi(Jb)$ by using gauge terms. To sum up, modulo ``canonical" gauge, the second order terms of $K$ can be only chosen as $\triangle J+\mathcal{R}$, and the gauge term is also unique, that is $\overline{X}=\overline{X_{1}}-\frac{1}{2}\overline{X_{2}}$.

Next, we consider deformation of $g$. We notice that as in the Ricci Flow $-2Ric+L_{\overline{X}}g$ gives a good symbol. And from the discussion above, the gauge terms are already chosen, so what we can do is just to find some ``canonical" second order tensors whose symbols compensate each other both in $g$ and $J$. But if it happens, it just gives the first order tensors. So modulo canonical gauge, $h$ is also unique, that is $-2Ric$. So we finish the proof of Theorem \ref{ah}.
\qed

\textit{Proof of Theorem \ref{ak}:}
In \aK setting, 
we require $\eta=P$. From Remark \ref{eta}, we have $K^{skew}=P^{(0,2)+(2,0)}, h^{(1,1)}=-P^{(1,1)}J.$ And from the calculations in \cite{ST symp}, we obtain
\begin{\eq}
h^{(1,1)}&=&-2Ric^{(1,1)}+\frac{1}{2}B^{1}-B^{2}\\
K^{skew}&=&\triangle J+\mathcal{N}.
\end{\eq}

So we have the freedom to choose the symmetric part of $K$. For the first order terms, from Lemma \ref{akt}, we see there is no such $(0,2)+(2,0)$ and symmetric tensors in \aK setting. For the second order terms, first we notice that the canonical gauge vanishes, and from Lemma \ref{cur} notice that $(\rho' J)^{sym}$ is $(1,1)$, we see the only candidate is $\mathcal{R}$. So we just need to investigate it from the consideration of symbol. Comparing to \aH condition, we have a extra condition $d\omega=0$. So we need to check what new symbol identities we can obtain from this condition. Notice that
\begin{\eq}
d\omega =0\Leftrightarrow DJ(i,j,k)+DJ(j,k,i)+DJ(k,i,j)=0.
\end{\eq}

Considering symbol, we obtain
\begin{\eq}
\sigma(DJ)(i,j,k)+\sigma(DJ)(j,k,i)+\sigma(DJ)(k,i,j)=0.
\end{\eq}

By direct calculation, we have
\begin{\eq}
\sigma(DJ)(i,j,k)&=&\xi(i)K(j,k)+\frac{1}{2}(h(Jj,k)\xi(i)+h(Jk,j)\xi(i)\\
&&+h(i,k)\xi(Jj)+h(i,Jk)\xi(j)-h(i,Jj)\xi(k)-h(i,j)\xi(Jk)).
\end{\eq}

So
\begin{\eq}
0&=&\xi(i)K(j,k)+\frac{1}{2}(h(Jj,k)\xi(i)+h(Jk,j)\xi(i)\\
&&+h(i,k)\xi(Jj)+h(i,Jk)\xi(j)-h(i,Jj)\xi(k)-h(i,j)\xi(Jk))\\
&&+\xi(j)K(k,i)+\frac{1}{2}(h(Jk,i)\xi(j)+h(Ji,k)\xi(j)\\
&&+h(j,i)\xi(Jk)+h(j,Ji)\xi(k)-h(j,Jk)\xi(i)-h(j,k)\xi(Ji))\\
&&+\xi(k)K(i,j)+\frac{1}{2}(h(Ji,j)\xi(k)+h(Jj,i)\xi(k)\\
&&+h(k,j)\xi(Ji)+h(k,Jj)\xi(i)-h(k,Ji)\xi(j)-h(k,i)\xi(Jj)).
\end{\eq}

By simplified the identity above, we have
\begin{\eq}
\xi(i)K(j,k)+\xi(j)K(k,i)+\xi(k)K(i,j)+h(k,Jj)\xi(i)+h(i,Jk)\xi(j)+h(j,Ji)\xi(k)=0
\end{\eq}

To consider the symbol of the second order $2$-tensors, we just need to take tensor product with $\xi$ and take trace or complex trace. Since we consider the symbol of $K$, we require it is of odd type. We have
\begin{\eq}
&&\xi(l)\xi(i)K(j,k)+\xi(l)\xi(j)K(k,i)+\xi(l)\xi(k)K(i,j)\\
&&+h(k,Jj)\xi(i)\xi(l)+h(i,Jk)\xi(j)\xi(l)+h(j,Ji)\xi(k)\xi(l)=0.
\end{\eq}

Considering the symmetries, we have the following cases to take trace or complex trace.

Let $l=i,j=a,k=b$, we obtain
\begin{eqnarray*}
|\xi|^{2}K(a,b)+\xi(a)K(b,\xi)+\xi(b)K(\xi,a)+|\xi|^{2}h(Ja,b)+\xi(a)h(\xi,Jb)+\xi(b)h(a,J\xi)=0.
\end{eqnarray*}

Let $l=Ji,j=Ja,k=b$, we obtain
\begin{eqnarray*}
-\xi(Ja)K(b,J\xi)+\xi(b)K(\xi,a)+h(\xi,b)\xi(Ja)+h(Ja,\xi)\xi(b)=0.
\end{eqnarray*}

Let $i=j,k=a,l=b$, we obtain
\begin{eqnarray*}
\xi(b)K(\xi,a)+\xi(b)K(a,\xi)+h(a,J\xi)\xi(b)+h(\xi,Ja)\xi(b)=0.
\end{eqnarray*}

Let $i=Jj,k=Ja,l=b$, we obtain
\begin{eqnarray*}
K(\xi,a)\xi(b)-K(a,\xi)\xi(b)+h(Ja,\xi)\xi(b)-h(J\xi,a)\xi(b)-trh\xi(Ja)\xi(b).
\end{eqnarray*}

Recall Lemma \ref{symbol}, we see that, in $\sigma(\triangle J)$ and $\sigma(\mathcal{R})$, the terms involving $|\xi|^{2}$ are in fact, $|\xi|^{2}(hJ)^{(0,2)+(2,0)}$. To cancel this term, the only possibly useful identity is the first one, but when consider the $(0,2)+(2,0)$ part, it gives nothing. So for the consideration of symbol, the choice is unique. So we finish the proof of Theorem \ref{ak}.
\qed

\textit{Proof of Theorem \ref{h}:}
In Hermitian case, we require $\delta J=K=0$.
Consider the second order terms in $h$, it should satisfy $h^{(0.2)+(2,0)}=0$ and $\sigma (h)$ itself is $Id$, not modulo gauge. From our classification results, the candidates are
\begin{\eq}
Ric(X,Y)+Ric(JX,JX),\quad Rm(JX,Y,i,Ji)-Rm(X,JY,i,Ji),\\
\quad D^{2}J(JX,i,Y,i)+D^{2}J(JY,i,X,i)-D^{2}J(X,i,JY,i)-D^{2}J(Y,i,JX,i).
\end{\eq}

Notice that $J$ is fixed, so it is always integrable. So the deformation condition is just that $h$ is $(1,1)$. From Lemma \ref{orisym}, the symbols of the candidates are computed as follows.
\begin{\eq}
&&\sigma(Ric(X,Y)+Ric(JX,JY))_{ab}\\
&=&\frac{1}{2}(\xi(a)h(b,\xi)+\xi(b)h(a,\xi)-\xi(a)\xi(b)trh-|\xi|^{2}h(a,b)\\
&&+\xi(Ja)h(Jb,\xi)+\xi(Jb)h(Ja,\xi)-\xi(Ja)\xi(Jb)trh-|\xi|^{2}h(Ja,Jb))\\
&=&\frac{1}{2}(\xi(a)h(b,\xi)+\xi(b)h(a,\xi)+\xi(Ja)h(Jb,\xi)+\xi(Jb)h(Ja,\xi)-\xi(Ja)\xi(Jb)trh-\xi(a)\xi(b)trh)\\
&&-|\xi|^{2}h(a,b).\\
\\
&&\sigma(Ric(JX,Y,i,Ji)-Ric(X,JY,i,Ji))_{ab}\\
&=&\xi(Ja)h(b,J\xi)-\xi(b)h(Ja,J\xi)-\xi(a)h(Jb,J\xi)+\xi(Jb)h(a,J\xi)\\
&=&\xi(Ja)h(b,J\xi)-\xi(b)h(a,\xi)-\xi(a)h(b,\xi)+\xi(Jb)h(a,J\xi).\\
\\
&&\sigma(D^{2}J(JX,i,Y,i)+D^{2}J(JY,i,X,i)-D^{2}J(X,i,JY,i)-D^{2}J(Y,i,JX,i))_{ab}\\
&=&\xi(Ja)\xi(Jb)trh+\xi(a)\xi(b)trh+\xi(Ja)K(b,\xi)+\xi(Jb)K(a,\xi)-\xi(a)K(Jb,\xi)-\xi(b)K(Ja,\xi)\\
&&+\xi(Ja)h(J\xi,b)+\xi(Jb)h(J\xi,a)-\xi(a)h(J\xi,Jb)-\xi(b)h(J\xi,Ja)\\
&=&\xi(Ja)\xi(Jb)trh+\xi(a)\xi(b)trh+\xi(Ja)h(J\xi,b)+\xi(Jb)h(J\xi,a)-\xi(a)h(\xi,b)-\xi(b)h(\xi,a).
\end{\eq}

So we see that the ``good" second order term is
\begin{\eq}
-2Ric^{(1,1)}-2(D^{2}J(J\cdot,i,\cdot,i))^{sym,(1,1)}.
\end{\eq}

Since the ``good" choice of symbol term is unique, from the calculation in \cite{ST her}, we see in fact the symbol term above can be given from $\frac{\partial }{\partial t}\omega=S$. The expression of $\mathcal{Q}$ is from Lemma \ref{eta} and Lemma \ref{tenh}. So we finish the proof of Theorem \ref{h}.
\qed
\begin{rem}\label{hergauge}
From \cite{ST gkg}, we see in Hermitian setting, modulo first order terms, $\triangle J+\mathcal{R}$ is just $L_{\theta^{\sharp}}J$. And in fact the above desired second order terms is the $(1,1)$ part of $-2Ric-L_{\theta^{\sharp}}g$, which coincides with our result in \aH setting.
\end{rem}
\section{Appendix: Hermitian Connection}
In this appendix, we review some basic results about Hermitian connection. For further study, one may refer \cite{Ga}.

Let $(g,J,\omega)$ be an almost Hermitian structure. Let $D$ be Levi-Civita connection and $\bigtriangledown$ be a linear connection. Let $\bigtriangledown=D+A$, i.e., $g(\bigtriangledown_{X}Y,Z)=g(D_{X}Y,Z)+A(X,Y,Z)$, where $A$ is a $3$-tensor.
\begin{lem}\label{connection}
\begin{\eq}
\bigtriangledown g=0&\Leftrightarrow& A(X,Y,Z)+A(X,Z,Y)=0.\\
\bigtriangledown J=0&\Leftrightarrow& A(X,JY,Z)+A(X,Y,JZ)+DJ(X,Y,Z)=0.
\end{\eq}
\end{lem}
\pf
\begin{\eq}
\nabla g(X,Y,Z)&=&Xg(Y,Z)-g(\nabla_{X}Y,Z)-g(Y,\nabla_{X}Z)\\
&=&Xg(Y,Z)-g(D_{X}Y,Z)-A(X,Y,Z)-g(Y,D_{X}Z)-A(X,Z,Y)\\
&=&-A(X,Y,Z)-A(X,Z,Y).\\
g(\nabla J(X,Y),Z)&=&g(\nabla_{X}(JY)-J\nabla_{X}Y,Z)\\
&=&g(D_{X}(JY),Z)+A(X,JY,Z)-g(JD_{X}Y,Z)+A(X,Y,JZ)\\
&=&DJ(X,Y,Z)+A(X,JY,Z)+A(X,Y,JZ).
\end{\eq}
\qed

If $\bigtriangledown g=\bigtriangledown J=0$, then we say that $\bigtriangledown$ is an Hermitian connection. In general, Hermitian connection is not unique. Naturally, we assume $A$ is defined from $(g,J,\omega)$. Since connection is of first order, we require $A$ is of first order. And since $D$ is of even type, we require $A$ is of even type. In a word, we assume $A=J*DJ$.
\begin{lem}
Let $\bigtriangledown$ be an Hermitian connection.

In almost Hermitian setting,
\begin{\eq}
A=\frac{1}{2}DJ(X,JY,Z)+\frac{t}{4}(DJ(JY,Z,X)+DJ(JZ,X,Y)-DJ(Y,Z,JX)-DJ(Z,X,JY)),
\end{\eq}
In Hermitian setting,
\begin{\eq}
A=\frac{1}{2}DJ(X,JY,Z)-\frac{t}{2}(DJ(Y,Z,JX)+DJ(Z,X,JY)).
\end{\eq}

In \aK setting,
\begin{\eq}
A=\frac{1}{2}DJ(X,JY,Z).
\end{\eq}
\end{lem}
\pf

Suppose $(g,J,\omega)$ is an almost Hermitian structure. From our assumption and Lemma \ref{dj}, we have
\begin{\eq}
A(X,Y,Z)&=&a_{1}DJ(X,Y,JZ)+a_{2}DJ(Y,Z,JX)+a_{3}DJ(Z,X,JY)\\
&&+a_{4}DJ(JX,Y,Z)+a_{5}DJ(JY,Z,X)+a_{6}DJ(JZ,X,Y).
\end{\eq}

From Lemma \ref{connection}, $\nabla$ is Hermitian if and only if
\begin{\eq}
0&=&a_{1}DJ(X,Y,JZ)+a_{2}DJ(Y,Z,JX)+a_{3}DJ(Z,X,JY)\\
&&+a_{4}DJ(JX,Y,Z)+a_{5}DJ(JY,Z,X)+a_{6}DJ(JZ,X,Y)\\
&&+a_{1}DJ(X,Z,JY)+a_{2}DJ(Z,Y,JX)+a_{3}DJ(Y,X,JZ)\\
&&+a_{4}DJ(JX,Z,Y)+a_{5}DJ(JZ,Y,X)+a_{6}DJ(JY,X,Z).\\
-DJ(X,Y,Z)&=&a_{1}DJ(X,JY,JZ)+a_{2}DJ(JY,Z,JX)-a_{3}DJ(Z,X,Y)\\
&&+a_{4}DJ(JX,JY,Z)-a_{5}DJ(Y,Z,X)+a_{6}DJ(JZ,X,JY)\\
&&-a_{1}DJ(X,Y,Z)+a_{2}DJ(Y,JZ,JX)+a_{3}DJ(JZ,X,JY)\\
&&+a_{4}DJ(JX,Y,JZ)+a_{5}DJ(JY,JZ,X)-a_{6}DJ(Z,X,Y).
\end{\eq}

To simplify the above equations, we have
\begin{\eq}
0&=&(a_{2}-a_{3})(DJ(Y,Z,JX)-DJ(Z,X,JY))+(a_{5}-a_{6})(DJ(JY,Z,X)-DJ(JZ,X,Y)).\\
0&=&(1-2a_{1})DJ(X,Y,Z)+2a_{4}DJ(JX,JY,Z)\\
&&+(a_{2}+a_{5})(DJ(JY,Z,JX)-DJ(Y,Z,X))+(a_{3}+a_{6})(DJ(JZ,X,JY)-DJ(Z,X,Y)).
\end{\eq}

Therefore, in \aH setting, $a_{1}=\frac{1}{2}$, $a_{4}=0$, $a_{2}=a_{3}=-a_{5}=-a_{6}=-\frac{t}{4}$.

In Hermitian setting, from Lemma \ref{hc},
\begin{\eq}
A(X,Y,Z)=(a_{1}-a_{4})DJ(X,Y,JZ)+(a_{2}-a_{5})DJ(Y,Z,JX)+(a_{3}-a_{6})DJ(Z,X,JY),
\end{\eq}

and the equations are
\begin{\eq}
0&=&(a_{2}-a_{5}-a_{3}+a_{6})(DJ(Y,Z,JX)-DJ(Z,X,JY))\\
0&=&(1-2a_{1}+2a_{4})DJ(X,Y,Z).
\end{\eq}

Therefore, $a_{1}-a_{4}=\frac{1}{2}$, $a_{2}-a_{5}=a_{3}-a_{6}=-\frac{t}{2}$.

In \aK setting, from Lemma \ref{akc},
\begin{\eq}
A(X,Y,Z)=(a_{1}+a_{4})DJ(X,Y,JZ)+(a_{2}+a_{5})DJ(Y,Z,JX)+(a_{3}+a_{6})DJ(Z,X,JY),
\end{\eq}

and the equations are
\begin{\eq}
0&=&(a_{2}+a_{5}-a_{3}-a_{6})(DJ(Y,Z,JX)-DJ(Z,X,JY))\\
0&=&(1-2a_{1}-2a_{4})DJ(X,Y,Z)-(2a_{2}+2a_{5})DJ(Y,Z,X)-(2a_{3}+2a_{6})DJ(Z,X,Y).
\end{\eq}

Therefore, $a_{2}+a_{5}=a_{3}+a_{6}$, $a_{1}+a_{4}=a_{2}+a_{5}+\frac{1}{2}$. Then
\begin{\eq}
A(X,Y,Z)=\frac{1}{2}DJ(X,Y,JZ).
\end{\eq}
\qed


\begin{thebibliography}{99}
\bibitem{B} Boling, J. $Homogeneous~ solutions~ of~ pluriclosed~ flow~ on~ closed~ complex~ surfaces$, arXiv preprint arXiv:1404.7106, 2014.
\bibitem{Bryant} http://mathoverflow.net/questions/72906/how-to-deduce-this-equation-for-a-4-dim-almost-kahler-manifold
\bibitem{Cao} Cao, H.D. $Deformation~ of~ K\ddot{a}hler~ metrics~ to~ K\ddot{a}hler-Einstein~ metrics~ on~ compact~ K\ddot{a}hler~ manifolds$, Invent. Math. 81 (1985), no. 2, 359--372.
\bibitem{Dai} Dai, S. $A~ Curvature~ flow~ unifying~ symplectic~ curvature~ flow~ and~ pluriclosed~ flow$, Pacific Journal of Mathematics, Vol. 277 (2015), No. 2, 287¨C311.
\bibitem{E} Enrietti, N. $Static~ SKT~ metrics~ on~ Lie~ groups$, Manuscripta Mathematica, 2013: 1-15.
\bibitem{EFV} Enrietti, N., Fino, A., Vezzoni, L. $The~ pluriclosed~ flow~ on~ nilmanifolds~ and~ Tamed~ symplectic~ forms$, arXiv preprint arXiv:1210.4816, 2012.
\bibitem{F} Fern\'{a}ndez-Culma, E. $Soliton~ almost~ K\ddot{a}hler~ structures~ on~ 6-dimensional~ nilmanifolds~ for~ the~
symplectic~ \\curvature~ flow$, arXiv preprint arXiv:1303.5461, 2013.
\bibitem{FV} Fino, A., Vezzoni, L. $Special~ Hermitian~ metrics~ on~ compact~ solvmanifolds$, Journal of Geometry and Physics.
\bibitem{Ga} Gauduchon, P. $Hermitian~ connections~ and~ Dirac~ operators$, Bollettino della Unione Matematica Italiana-B, 1997 (2): 257-288.
\bibitem{Gill} Gill, M. $Convergence~ of~ the~ parabolic~ complex~ Monge-Amp\grave{e}re~ equation~ on~ compact~ Hermitian~ manifolds$, Comm. Anal. Geom 19(2011), no.2, 277-303.
\bibitem{Ham} Hamilton, R. $Three-manifolds~ with~ positive~ Ricci~ curvature$, Journal of Differential Geometry 17 (1982), no. 2, 255--306.
\bibitem{Lau} Lauret, J. $Curvature~ flows~ for~ almost-Hermitian~ Lie~ groups$, Trans. Amer. Math. Soc. 367 (2015), 7453-7480.
\bibitem{LW} Lauret, J., Will, C. $On~ the~ symplectic~ curvature~ flow~ for~ locally~ homogeneous~ manifolds$ arXiv preprint arXiv:1405.6065, 2014.
\bibitem{LeWang} L\^{e}, H.V., Wang, G.F. $Anti-complexified~ Ricci~ flow~ on~ compact~ symplectic~ manifolds$, J. Reign Angew. Math. 530(2011), 17-31.
\bibitem{NN} Newlander, A., Nirenberg, L. $Complex~ analytic~ coordinates~ in~ almost~ complex~ manifolds$, The Annals of Mathematics, 1957, 65(3): 391-404.
\bibitem{Pe1} Perelman, G. $The~ entropy~ formula~ for~ the~ Ricci~ flow~ and~ its~ geometric~ applications$, arXiv preprint math/0211159, 2002.
\bibitem{Pe2} Perelman, G. $Ricci~ flow~ with~ surgery~ on~ three-manifolds$, arXiv preprint math/0303109, 2003.
\bibitem{Pe3} Perelman, G. $Finite~ extinction~ time~ for~ the~ solutions~ to~ the~ Ricci~ flow~ on~ certain~ three-manifolds$, arXiv preprint math/0307245, 2003.
\bibitem{Po} Pook, J. $Homogeneous~ and~ locally~ homogeneous~ solutions~ to~ symplectic~ curvature~ flow$, arXiv preprint arXiv:1202.1427, 2012.
\bibitem{Sm} Smith, D. $Stability~ of~ the~ Almost~ Hermitian~ Curvature~ Flow$, arXiv preprint arXiv:1308.6214, 2013.
\bibitem{Streets pgs} Streets, J. $Pluriclosed~ flow~ on~ generalized~ K\ddot{a}hler~ manifolds~ with~ split~ tangent~ bundle$, to appear in Crelle's Journal.
\bibitem{Streets pbg} Streets, J. $Pluriclosed~ flow, ~Born-Infeld ~geometry,~ and~ rigidity~ results~ for ~ generalized~ K\ddot{a}hler~ manifolds$, arXiv preprint arXiv:1502.02584, 2015.
\bibitem{Streets gkrf} Streets, J. $Generalized~ K\ddot{a}hler-Ricci~ flow~ and~ the~ classification~ of~ nondegenerate~ generalized~ K\ddot{a}hler~ surfaces$, arXiv preprint arXiv:1601.02981, 2016.
\bibitem{ST pluri} Streets, J., Tian, G. $A~ parabolic~ flow~ of~ pluriclosed~ metrics$, Int. Math. Res. Notices (2010), (16): 3101-3133.
\bibitem{ST reg pluri} Streets, J., Tian, G. $Regularity~ results~ for~ pluriclosed~ flow$, to appear in Geometry \& Topology.
\bibitem{ST her} Streets, J., Tian, G. $Hermitian~ curvature~ flow$, Journal of the European Mathematical Society, 2011, 13(3): 601-634.
\bibitem{ST symp} Streets, J., Tian, G. $Symplectic~ curvature~ flow$, Journal f\"{u}r die reine und angewandte Mathematik (Crelles Journal), 2011.
\bibitem{ST gkg} Streets, J., Tian, G. $Generalized~ K\ddot{a}hler~ geometry~ and~ the~ pluriclosed~ flow$, Nuclear Physics B, 2012, 858(2): 366-376.
\bibitem{SW} Streets, J., Warren, M. $Evans-Krylov~ Estimates~ for~ a~ nonconvex~ Monge-Amp\grave{e}re~ equation$, arXiv preprint arXiv:1410.2911, 2014.
\bibitem{TW} Tossati, V., Weinkove, B. $On~ the~ evolution~ of~ a~ Hermitian~ metric~ by~ its~ Chern-Ricci~ form$, arXiv preprint arXiv:1201.0312, 2012.
\bibitem{V} Vezzoni, L. $On~ Hermitian~ curvature~ flow~ on~ almost~ complex~ manifolds$, Differential Geometry and its Applications, 2011, 29(5): 709-722.
\end{thebibliography}
\end{document}